\documentclass[amsfonts, 12pt]{article}

\usepackage{amsmath,amsthm,amsfonts,amssymb}

\newtheorem{theorem}{Theorem}[section]

\newtheorem{example}[theorem]{Example}

\newtheorem{lemma}[theorem]{Lemma}
\newtheorem{definition}[theorem]{Definition}
\newtheorem{remark}[theorem]{Remark}

\def\cB{\mathcal{B}}

\def\cE{\mathcal{E}}

\def\cH{\mathcal{H}}
\def\cF{\mathcal{F}}
\def\cL{\mathcal{L}}

\def\cU{\mathcal{U}}

\def\bC{\mathbb{C}}
\def\bD{\mathbb{D}}

\def\bR{\mathbb{R}}

\begin{document}

\title{Malliavin differentiability of solutions of SPDEs with L\'evy white noise}

\author{Raluca M. Balan\footnote{University of Ottawa, Department of Mathematics and Statistics,
585 King Edward Avenue, Ottawa, ON, K1N 6N5, Canada. E-mail
address: rbalan@uottawa.ca} \footnote{Research supported by a
grant from the Natural Sciences and Engineering Research Council
of Canada.}\and
Cheikh B. Ndongo\footnote{University of Ottawa, Department of Mathematics and Statistics,
585 King Edward Avenue, Ottawa, ON, K1N 6N5, Canada. E-mail
address: cndon072@uottawa.ca}
}

\date{May 9, 2016}
\maketitle

\begin{abstract}
\noindent In this article, we consider a stochastic partial differential equation (SPDE) driven by a L\'evy white noise, with Lipschitz multiplicative term $\sigma$. We prove that under some conditions, this equation has a unique random field solution. These conditions are verified by the stochastic heat and wave equations.
We introduce the basic elements of Malliavin calculus with respect to
the compensated Poisson random measure associated with the L\'evy white noise.
If $\sigma$ is affine, we prove that the solution is Malliavin differentiable and its Malliavin derivative satisfies a stochastic integral equation.
\end{abstract}

\vspace{10mm}


\vspace{5mm}

%60H15 s.p.d.e.s
%60G51 processes with independent increments and Levy processes

\pagebreak

\section{Introduction}

In this article, we consider the stochastic partial differential equation (SPDE):
\begin{equation}
\label{spde}
\cL u(t,x) =  \sigma(u(t,x))\dot{L}(t,x), \quad t\in [0,T],x \in \bR
\end{equation}
with some deterministic initial conditions, where $\cL$ is a second-order differential operator on $[0,T] \times \bR$, $\dot{L}$ denotes the formal derivative of the L\'evy white noise $L$ (defined below) and the function $\sigma:\bR \to \bR$ is  Lipschitz continuous.

A process $u=\{u(t,x);t \in [0,T],x \in \bR\}$ is called a (mild) {\em solution} of \eqref{spde} if $u$ is predictable and satisfies the following integral equation:
$$u(t,x)=w(t,x)+\int_0^t \int_{\bR}G(t-s,x-y)\sigma(u(s,y))L(ds,dy),$$
where $w$ is the solution of the deterministic equation $\cL u=0$ with the same initial conditions as \eqref{spde} and $G$ is the Green function of the operator $\cL$.

The study of SPDEs with Gaussian noise is a well-developed area of stochastic analysis, and the behaviour of random-field solutions of such equations is well-understood. We refer the reader to \cite{walsh86} for the original lecture notes which lead to the development of this area, and to \cite{DKMNX,K14} for some recent advances. In particular, the probability laws of these solutions can be analyzed using techniques from Malliavin calculus, as described in \cite{nualart06,sanz05}.

On the other hand, there is a large literature dedicated to the study of stochastic differential equations (SDE) with L\'evy noise, the monograph \cite{applebaum09} containing a comprehensive account on this topic. One can develop also a Malliavin calculus for L\'evy processes with finite variance, using an analogue of the Wiener chaos representation with respect to underlying Poisson random measure of the L\'evy process. This method was developed in \cite{BGJ87} with the same purpose of analyzing the probability law of the solution of an SDE driven by a finite variance L\'evy noise. More recently, Malliavin calculus for L\'evy processes with finite variance have been used in financial mathematics, the monograph \cite{DOP09} being a very readable introduction to this topic.

There are two approaches to SPDEs in the literature. One is the random field approach
which originates in John Walsh's lecture notes \cite{walsh86}. When using this approach, the solution is viewed as a real-valued process which is indexed by time and space. The other approach is the infinite-dimensional approach, due to Da Prato and Zabczyk \cite{DZ92}, according to which the solution is a process indexed by time only, which takes values in an infinite-dimensional Hilbert space. It is not always possible to compare the solutions obtained using the two approaches (see  \cite{dalang-quer11} for several results in this direction).
SPDEs with L\'evy noise were studied in the monograph \cite{PZ07}, using the infinite-dimensional approach. In the present article, we use the random field approach for examining an SPDE driven by the finite variance L\'evy noise introduced in \cite{B15}, with the goal of studying the Malliavin differentiability of the solution. As mentioned above, this study can be useful for analyzing the probability law of the solution. We postpone this problem for future work.

We begin by recalling from \cite{B15} the construction of the L\'evy white noise $L$ driving equation \eqref{spde}.
We consider a Poisson random measure (PRM) $N$ on the space $U=[0,T] \times \bR \times \bR_0$ of intensity $\mu=dt dx \nu(dz)$ defined on a complete probability space $(\Omega, \cF,P)$, where $\nu$ is a L\'evy measure on $\bR_0$, i.e. $\nu$ satisfies
$$\int_{\bR_0}(1 \wedge |z|^2)\nu(dz)<\infty.$$
Here $\bR_0=\bR \verb2\2\{0\}$. In addition, we assume that $\nu$ satisfies the following condition:
$$v:=\int_{\bR_0}z^2 \nu(dz)<\infty.$$

We denote by $\widehat{N}$ the compensated PRM defined by $\widehat{N}(A)=N(A)-\mu(A)$ for any $A \in \cU$ with $\mu(A)<\infty$, where $\cU$ is the class of Borel sets in $U$. We denote by $\cF_t$ the $\sigma$-field generated by $N([0,s] \times B \times \Gamma)$ for all $s \in [0,t],B \in \cB_b(\bR)$ and $\Gamma \in \cB_b(\bR_0)$. We denote by $\cB_b(\bR)$ the class of bounded Borel sets in $\bR$, and by $\cB_b(\bR_0)$ the class of Borel sets in $\bR_0$ which are bounded away from $0$.

A {\em L\'evy white noise} with intensity measure $\nu$ is a collection $L=\{L_t(B);t \in [0,T], B \in \cB_b(\bR)\}$ of zero-mean square-integrable random variables defined by
$$L_t(B)=\int_0^t \int_B z \widehat{N}(ds,dx,dz).$$
These variables have the following properties:\\
(i) $L_0(B)=0$ a.s. for all $B \in \cB_b(\bR)$;\\
(ii) $L_t(B_1), \ldots,L_t(B_k)$ are independent for any $t >0$ and for any disjoint sets $B_1, \ldots, B_k \in \cB_b(\bR)$;\\
(iii) for any $0<s\leq t$ and for any $B \in \cB_b(\bR)$, $L_t(B)-L_s(B)$ is independent of $\cF_s$ and has characteristic function
$$E(e^{iu(L_t(B)-L_s(B))})=\exp\left\{(t-s)|B|\int_{\bR_0} (e^{iuz}-1-iuz)\nu(dz)\right\}, \quad u \in \bR.$$
We denote by $\cF_t^L$ the $\sigma$-field generated by $L(s)$ for all $s \in [0,t]$. For any $h \in L^2([0,T] \times \bR)$, we define the stochastic integral of $h$ with respect to $L$:
$$L(h)=\int_0^T \int_{\bR}h(t,x)L(dt,dx)=\int_0^T \int_{\bR}\int_{\bR_0}h(t,x)z\widehat{N}(dt,dx,dz).$$

Using the same method as in It\^o's classical theory, this integral can be extended to
random integrands, i.e. to the class of predictable processes $X=\{X(t,x); t \in [0,T], x\in \bR\}$ such that $E \int_0^T \int_{\bR}|X(t,x)|^2 dxdt<\infty$. The integral has the following isometry property:
\begin{equation}
\label{isometry}
E\left|\int_0^T \int_{\bR}X(t,x)L(dt,dx)\right|^2=v E \int_0^T \int_{\bR}|X(t,x)|^2 dxdt.
\end{equation}
Recall that a process $X=\{X(t,x); t \geq 0, x\in \bR^d\}$ is {\em predictable} if it is measurable with respect to the predictable $\sigma$-field on $\bR_+ \times \bR$, i.e. the $\sigma$-field generated by processes of the form
$X(\omega,t,x)=Y(\omega)1_{(a,b]}(t) 1_{A}(x)$, where $0<a<b$, $Y$ is a bounded and $\cF_a^L$-measurable random variable and $A \in \cB_b(\bR)$.

This article is organized as follows. In Section \ref{section-Malliavin}, we introduce the basic elements of Malliavin calculus with respect to the compensated Poisson random measure $\widehat{N}$. In Section \ref{section-exist}, we prove that under a certain hypothesis, equation \eqref{spde} has a unique solution. This hypothesis is verified in the case of the wave and heat equations.
In Section \ref{section-Mall-sol}, we examine the Malliavin differentiability of the solution, in the case when the function $\sigma$ is affine. Finally, in Appendix \ref{section-app}, we include a version of Gronwall's lemma which is needed in the sequel.

\section{Malliavin calculus on the Poisson space}
\label{section-Malliavin}

In this section, we introduce the basic ingredients of Malliavin calculus with respect to the $\widehat{N}$, following very closely the approach presented in Chapters 10-12 of \cite{DOP09}. The difference compared to \cite{DOP09} is that our parameter space $U$ has variables $(t,x,z)$ instead of $(t,z)$. For the sake of brevity, we do not include the proofs of the results presented in this section. These proofs can be found in Chapter 6 of the doctoral thesis \cite{ndongo16} of the second author.

We let $\cH=L^2(U,\cU,\mu)$ and $\cH^{\otimes n}=L^2(U^n, \cU^n, \mu^n)$. We denote by $\cH^{\odot}$ the set of all symmetric functions $f \in \cH^{\otimes n}$.
We denote by $\cH_{\bC}, \cH_{\bC}^{\otimes n},\cH_{\bC}^{\odot n}$ the analogous spaces of $\bC$-valued functions.

Let $S_n=\{(u_1, \ldots,u_n) \in U^n;u_i=(t_i,x_i,z_i) \ \mbox{with} \ t_1<\ldots<t_n\}$. For any measurable function $f:S_n \to \bR$ with
$$\|f\|_{L^2(S_n)}:=\int_{S_n}|f(u_1,\ldots,u_n)|^2 d\mu^n(u_1, \ldots,u_n)<\infty,$$
we define the {\em $n$-fold iterated integral} of $f$ with respect to $\widehat{N}$ by
$$J_n(f)=\int_0^T \int_{\bR} \int_{\bR_0} \left(\int_{0}^{t_n-}\int_{\bR} \int_{\bR_0}  \ldots
 \left(\int_{0}^{t_2-}\int_{\bR} \int_{\bR_0} f(u_1,\ldots,u_n) \widehat{N}(du_1) \right)
 \ldots \right)\widehat{N}(du_n)$$
where $u_i=(t_i,x_i,z_i)$. Then $E[J_n(f)J_m(g)]=0$ for all $n\not=m$ and $E|J_n(f)|^2=n! \|f\|_{L^2(S_n)}^2$.

For any $f \in \cH^{\odot n}$ we defined the {\em multiple integral} of $f$ with respect to $\widehat{N}$ by
$I_n(f)=n! J_n(f)$. It follows that $E[I_n(f)I_m(g)]=0$ for all $n\not=m$ and
$$E|I_n(f)|^2=n! \|f\|_{\cH^{\otimes n}}^{2} \quad \mbox{for all} \quad f \in \cH^{\odot n}.$$
If $f\in \cH_{\bC}^{\odot n}$ with $f=g+ih$, we define $I_n(f)=I_n(g)+iI_n(h)$.

%can be represented as
%$$F=E(F)+\int_{0}^{T}\int_{\bR}\int_{\bR_0}\psi(t,x,z)\widehat{N}(dt,dx,dz),$$
%for some predictable $\sigma$-field $\psi \in \cH$. Based on this %representation, it follows that for any random variable $F \in L_{\bC}^2(\Omega,\cF_T^L,P)$

Let $L_{\bC}^{2}(\Omega)$ be the set of $\bC$-valued square-integrable random variables defined on $(\Omega,\cF,P)$.
By Theorem 7 of \cite{balan-ndongo15}, any $\cF_T^L$-measurable random variable $F \in L_{\bC}^{2}(\Omega)$ admits the chaos expansion
$$F=\sum_{n \geq 0}I_n(f_n) \quad \mbox{in} \quad L^2_{\bC}(\Omega),$$
where $f_n \in \cH^{\odot n}$ for all $n \geq 1$ and $f_0=E(F)$.

The chaos expansion plays a crucial role in developing the Malliavin calculus with respect to $\widehat{N}$. In particular, the Skorohod integrals with respect to $\widehat{N}$ and $L$ are defined as follows.

\begin{definition}
{\rm a) Let $X=\{X(u);u\in U\}$ be a square-integrable process such that $X(u)$ is $\cF_T^L$-measurable for any $u \in U$. For each $u \in U$ let
$X(u)=\sum_{n \geq 0}I_n(f_n(\cdot,u))$
be the chaos expansion of $X(u)$, with $f_n(\cdot,u) \in \cH^{\odot n}$. We denote by $\widetilde{f}_n(u_1, \ldots,u_n,u)$ the symmetrization of $f_n$ with respect to all $n+1$ variables.
We say that $X$ is {\em Skorohod integrable} with respect to $\widehat{N}$ (and we write $X \in {\rm Dom}(\delta)$) if
$$\sum_{n \geq 0}E|I_{n+1}(\widetilde{f}_n)|^2=\sum_{n \geq 1}(n+1)!\|\widetilde{f}_n\|_{\cH^{\odot(n+1)}}^{2}<\infty.$$
In this case, we define the {\em Skorohod integral} of $X$ with respect to $\widehat{N}$ by
$$\delta(X)=\int_{0}^{T}\int_{\bR}\int_{\bR_0}X(t,x,z)\widehat{N}(\delta t,\delta x, \delta z):=\sum_{n \geq 1} I_{n+1}(\widetilde{f}_{n}).$$

b) Let $Y=\{Y(t,x);t \in [0,T],x \in \bR\}$ be a square-integrable process such that $Y(t,x)$ is $\cF_T^L$-measurable for any $t \in [0,T]$ and $x \in \bR$. We say that $Y$ is {\em Skorohod integrable} with respect to $L$ (and we write $Y \in {\rm Dom}(\delta^L)$) if the process $\{Y(t,x)z; (t,x,z)\in U\}$ is Skorohod integrable with respect to $\widehat{N}$. In this case, we define the {\em Skorohod integral} of $Y$ with respect to $L$ by
$$\delta^L(Y)=\int_0^T \int_{\bR}Y(t,x)L(\delta t,\delta x):=\int_0^T \int_{\bR}\int_{\bR_0}Y(t,x)z\widehat{N}(\delta t,\delta x,\delta z).$$
}
\end{definition}

The following result shows that the Shorohod integral can be viewed as an extension of the It\^o integral.

\begin{theorem}
\label{Ito-Skorohod}
a) If $X=\{X(u);(u)\in U\}$ is a predictable process such that
 $E\|X\|_{U}^{2}<\infty$, then $X$ is Skorohod integrable with respect to $\widehat{N}$ and
$$\int_0^T\int_{\bR}\int_{\bR_0}X(t,x,z)\widehat{N}(\delta t,\delta x,\delta z)=\int_0^T \int_{\bR}\int_{\bR_0}X(t,x,z)\widehat{N}(dt,dx,dz).$$

b) If $Y=\{Y(t,x);t \in [0,T],x\in \bR\}$ is a predictable process such that $E\int_0^T \int_{\bR}|Y(t,x)|^2dxdt<\infty$, then $Y$ is Skorohod integrable with respect to $L$ and
$$\int_0^T \int_{\bR}Y(t,x)L(\delta t,\delta x)=\int_0^T \int_{\bR}Y(t,x)L(dt,dx).$$
\end{theorem}

We now introduce the definition of the Malliavin derivative.
 %with respect to $\widehat{N}$.

\begin{definition}
{\rm Let $F \in L^2(\Omega)$ be an $\cF_T^L$-measurable random variable with
the chaos expansion
$F=\sum_{n \geq 0}I_n(f_n)$ with $f_n \in \cH^{\odot n}$.
We say that $F$ is {\em Malliavin differentiable} with respect to $\widehat{N}$ if
$$\sum_{n \geq 1}n n! \|f_n\|_{\cH^{\otimes n}}^{2}<\infty.$$
In this case, we define the {\em Malliavin derivative} of $F$ with respect to $\widehat{N}$ by
$$D_{u}F=\sum_{n \geq 1}n I_{n-1}(f_n(\cdot,u)), \quad u \in U.$$
We denote by $\bD^{1,2}$ the space of Malliavin differentiable random variables with respect to $\widehat{N}$.}
\end{definition}

Note that $E\|D F\|_{\cH}^{2}=\sum_{n \geq 1}n n! \|f_n\|_{\cH^{\otimes n}}^{2}<\infty$.
%For any $p \geq 2$, we let
%$$\|F\|_{\bD^{1,p}}^{p}:=E|F|^p+E\|\bD_{\cdot}F\|_{\cH}^{p}$$
%and we denote by $\bD^{1,p}$ the space of all $\cF_{T}^L$-measurable random %variables $F$ such that $\|F\|_{\bD^{1,p}}<\infty$.

\begin{theorem}[Closability of Malliavin derivative]
\label{closability-D}
Let $(F_n)_{n \geq 1} \subset \bD^{1,2}$ and $F \in L^2(\Omega)$ such that $F_n \to F$ in $L^2(\Omega)$ and $(DF_n)_{n \geq 1}$ converges in $L^2(\Omega;\cH)$. Then $F \in \bD^{1,2}$ and $DF_n \to DF$ in $L^2(\Omega;\cH)$.
\end{theorem}

Typical examples of Malliavin differentiable random variables are  exponentials of stochastic integrals: for any $h \in L^2([0,T] \times \bR)$, $$D_{t,x,z}(e^{L(h)})=e^{L(h)}(e^{h(t,x)z}-1).$$ Moreover, the set $\bD_{\cE}^{1,2}$ of linear combinations of random variables of the form $e^{L(h)}$ with $h \in L^2([0,T] \times \bR)$ is dense in $\bD^{1,2}$.

The following result shows that the Malliavin derivative is a difference operator with respect to $\widehat{N}$, not a differential operator.

\begin{theorem}[Chain Rule]
\label{chain-rule}
For any $F \in \bD^{1,2}$ and any continuous function $g:\bR \to \bR$ such that $g(F) \in L^2(\Omega)$ and $g(F+DF)-g(F) \in L^2(\Omega;\cH)$,
$g(F) \in \bD^{1,2}$ and
$$Dg(F)=g(F+D F)-g(F) \quad \mbox{in} \quad L^2(\Omega;\cH).$$
\end{theorem}

Similarly to the Gaussian case, we have the following results.

\begin{theorem}[Duality Formula] If $F \in \bD^{1,2}$ and $X \in {\rm Dom}(\delta)$, then
$$E\left[F \int_0^T \int_{\bR}\int_{\bR_0}X(t,x,z)\widehat{N}(\delta t,\delta x,\delta z) \right]=E \left[\int_0^T \int_{\bR}\int_{\bR_0}X(t,x,z)D_{t,x,z}F \nu(dz)dxdt \right].$$
\end{theorem}

\begin{theorem}[Fundamental Theorem of Calculus]
Let $X=\{X(s,y,\zeta);s \in [0,T],y \in \bR,\zeta \in \bR_0\}$ be a process which satisfies the following conditions:\\
(i) $X(s,y,\zeta) \in \bD^{1,2}$ for any $(s,y,\zeta) \in U$;\\
(ii) $E\int_0^T \int_{\bR}\int_{\bR_0}|X(s,y,\zeta)|^2 \nu(dz)dyds<\infty$;\\
(iii) $\{D_{t,x,z}X(s,y,\zeta); (s,y,\zeta) \in U\} \in {\rm Dom}(\delta)$ for any $(t,x,z)\in U$;\\
(iv) $\{\delta(D_{t,x,z}X); (t,x,z) \in U \} \in L^2(\Omega;\cH)$.\\
Then $X \in {\rm Dom}(\delta)$, $\delta(X) \in \bD^{1,2}$ and $D[\delta(X)]=X+\delta(DX)$, i.e.
\begin{eqnarray*}
\lefteqn{D_{t,x,z}\left(\int_0^T \int_{\bR} \int_{\bR_0}X(s,y,\zeta)\widehat{N}(\delta s,\delta y,\delta \zeta) \right)= X(t,x,z)+} \\
& & \int_0^T \int_{\bR}\int_{\bR_0}D_{t,x,z}X(s,y,\zeta)\widehat{N}(\delta s,\delta y,\delta \zeta) \quad \mbox{in} \quad L^2(\Omega; \cH).
\end{eqnarray*}
\end{theorem}

As an immediate consequence of the previous theorem, we obtain the following result.
\begin{theorem}
\label{FTC}
Let $Y=\{Y(s,y);s \in [0,T],y \in \bR\}$ be a process which satisfies the following conditions:\\
(i) $Y(s,y) \in \bD^{1,2}$ for all $s \in [0,T]$ and $y \in \bR$;\\
(ii) $E \int_0^T \int_{\bR}|Y(s,y)|^2 dyds<\infty$;\\
(iii) $\{D_{t,x,z}Y(s,y);s \in [0,T],y \in \bR\} \in {\rm Dom}(\delta^L)$ for any $(t,x,z) \in U$;\\
(iv) $E \int_0^T \int_{\bR} \int_{\bR_0} \left|\int_0^T \int_{\bR}D_{t,x,z}Y(s,y)L(\delta s,\delta y) \right|^2 \nu(dz)dxdt<\infty$.\\
Then $Y \in {\rm Dom}(\delta^L)$, $\delta^L(Y) \in \bD^{1,2}$ and the following relation holds in $L^2(\Omega; \cH)$:
$$D_{t,x,z}(\delta^L(Y))= Y(t,x)z+ \int_0^T \int_{\bR}\int_{\bR_0}D_{t,x,z}Y(s,y)L(\delta s,\delta y).$$
\end{theorem}

\section{Existence of Solution}
\label{section-exist}

In this section, we show that equation \eqref{spde} has a unique solution.

We recall that $w$ is the solution of the homogeneous equation ${\cal L}u=0$ with the same initial conditions as \eqref{spde}, and $G$ is the Green function of the operator ${\cal L}$ on $\bR_{+} \times \bR$. We assume that for any $t \in [0,T]$,
$G(t,\cdot) \in L^1(\bR)$ and we denote by $\cF G(t,\cdot)$ its Fourier transform:
$$\cF G(t,\cdot)(\xi)=\int_{\bR}e^{-i \xi x}G(t,x)dx.$$
We suppose that the following hypotheses holds:

\vspace{2mm}

{\em Hypothesis H1.} $w$ is continuous and uniformly bounded on $[0,T] \times \bR$.

\vspace{2mm}

{\em Hypothesis H2.} %The Green function $G$ of the operator $\cL$ satisfies:\\
a) $\int_0^T \int_{\bR}G^2(t,x)dxdt<\infty$;\\
b) the function $t \mapsto \cF G(t,\cdot)(\xi)$ is continuous on $[0,T]$, for any $\xi \in \bR^d$;\\
c) there exists $\varepsilon>0$ and a non-negative function $k_t(\cdot)$ such that
$$|\cF G(t+h,\cdot)(\xi)-\cF G(t,\cdot)(\xi)|\leq k_t(\xi)$$ for any $t \in [0,T]$ and $h \in [0,\varepsilon]$, and $\int_0^T \int_{\bR}|k_t(\xi)|^2 d\xi dt<\infty$.

Since $\sigma$ is a Lipschitz continuous function, there exists a constant $C_{\sigma}>0$ such that for any $x,y \in \bR$,
\begin{equation}
\label{Lip1}
|\sigma(x)-\sigma(y)| \leq C_{\sigma}|x-y|.
\end{equation}
In particular, for any $x \in \bR$,
\begin{equation}
\label{Lip2}
|\sigma(x)| \leq D_{\sigma}(1+|x|),
\end{equation}
where $D_{\sigma}=\max\{C_{\sigma},|\sigma(0)|\}$.

The following theorem is an extension of Theorem 1.1.(a) of \cite{balan-ndongo16} to an arbitrary operator $\cL$. The proof of this theorem (for the stochastic wave equation) was omitted from \cite{balan-ndongo16}. We include the proof here.

\begin{theorem}
Equation \eqref{spde} has a unique solution $u=\{u(t,x); t \in [0,T],x\in \bR^d\}$ which is $L^2(\Omega)$-continuous and satisfies
$$\sup_{(t,x) \in [0,T] \times \bR}E|u(t,x)|^2<\infty.$$
\end{theorem}

\noindent {\bf Proof:} {\em Existence.} We use the same argument as in the proof of Theorem 13 of \cite{dalang99}. We denote by $(u_n)_{n \geq 0}$ the sequence of Picard iterations defined by: $u_0(t,x)=w(t,x)$ and
\begin{equation}
\label{def-Picard}
u_{n+1}(t,x)= w(t,x)+\int_0^t \int_{\bR}G(t-s,x-y)\sigma(u_n(s,y))L(ds,dy), \quad n \geq 0.
\end{equation}
By induction on $n$, it can be proved that the following property holds:
\begin{equation}
\label{propertyP}
\tag{P}
\left\{
\begin{array}{rcl}
& (i)  &  \ \displaystyle u_n(t,x) \ \mbox{is well-defined for any} \ (t,x) \in [0,T] \times \bR,  \quad  \\[1ex]
& (ii)  &  \displaystyle K_n:=\sup_{(t,x) \in [0,T] \times \bR}E|u_n(t,x)|^{2}<\infty, \\ [1ex]
& (iii) & \displaystyle (t,x) \mapsto u_n(t,x) \ \mbox{is $L^2(\Omega)$-continuous on $[0,T] \times \bR$,} \\ [1ex]
& (iv) & \displaystyle u_n(t,x) \ \mbox{is $\cF_t$-measurable for any $t \in [0,T]$ and $x \in \bR$}.
\end{array} \right.
\end{equation}
(Hypotheses (H1) and (H2) are needed for the proof of property $(iii)$.)
From properties $(iii)$ and $(iv)$, it follows that $u_n$ has a predictable modification, denoted also by $u_n$. This modification is used in the definition \eqref{def-Picard} of $u_{n+1}(t,x)$. Using the isometry property \eqref{isometry} of the stochastic integral and \eqref{Lip1}, we have:
\begin{eqnarray*}
E|u_{n+1}(t,x)-u_n(t,x)|^2 &=& v E \int_0^t \int_{\bR}G^2(t-s,x-y) |\sigma(u_n(s,y))-\sigma(u_{n-1}(s,y))|^2 dy ds \\
& \leq & v C_{\sigma}^2 \int_{0}^t \int_{\bR}G^2(t-s,x-y)E|u_n(s,y)-u_{n-1}(s,y)|^2 dy ds \\
& \leq & v C_{\sigma}^2 \int_0^t M_n(s) \left(\int_{\bR}G^2(t-s,x-y)dy \right)ds,
\end{eqnarray*}
where $H_n(t)=\sup_{x \in \bR}E|u_n(t,x)-u_{n-1}(t,x)|^2$. For any $t \in [0,T]$, we denote
\begin{equation}
\label{def-J}
J(t)=\int_{\bR}G^2(t,x)dx.
\end{equation}
Taking the supremum over $x \in \bR$ in the previous inequality, we obtain that:
$$H_{n+1}(t) \leq v C_{\sigma} \int_{0}^t H_n(s)J(t-s)ds,$$
for any $t \in [0,T]$ and $n \geq 0$.
By applying Lemma 15 of \cite{dalang99} with $k_1=k_2=0$, we infer that
\begin{equation}
\label{sum-converges}
\sum_{n \geq 0}\sup_{t \in [0,T]}H_n(t)^{1/2}<\infty.
\end{equation}
This shows that the sequence $(u_n)_{n \geq 0}$ converges in $L^2(\Omega)$ to a random variable $u(t,x)$, uniformly in $[0,T] \times \bR$, i.e.
\begin{equation}
\label{un-conv-u}
\sup_{(t,x) \in [0,T] \times \bR}E|u_n(t,x)-u(t,x)|^2 \to 0.
\end{equation}
To see that $u$ is a solution of \eqref{spde}, we take the limit in $L^2(\Omega)$ as $n \to \infty$ in \eqref{def-Picard}. In particular, this argument shows that
\begin{equation}
\label{def-K}
K:=\sup_{n \geq 1}\sup_{(t,x) \in [0,T] \times \bR}E|u_n(t,x)|^2<\infty.
\end{equation}

{\em Uniqueness.} Let $H(t)=\sup_{x \in \bR}E|u(t,x)-u'(t,x)|^2$, where $u$ and $u'$ are two solutions of \eqref{spde}. A similar argument as above shows that
$$H(t)\leq v C_{\sigma}^2 \int_0^t H(s)J(t-s)ds,$$
for any $t \in [0,T]$. By Gronwall's lemma, $H(t)=0$ for all $t \in [0,T]$.
$\Box$

\begin{example}[wave equation]
{\rm If $\cL=\frac{\partial}{\partial^2 t}-\frac{\partial}{\partial x^2}$, then
$G(t,x)=\frac{1}{2}1_{\{|x| \leq t\}}$. Hypothesis (H2) holds since
$$\cF G(t,\cdot)(\xi)=\frac{\sin(t|\xi|)}{|\xi|}.$$
}
\end{example}

\begin{example}[heat equation]
{\rm If $\cL=\frac{\partial}{\partial t}-\frac{1}{2}\frac{\partial}{\partial x^2}$, then
$G(t,x)=(2 \pi t)^{-1/2}\exp\left(-\frac{|x|^2}{2t}\right)$. Hypothesis (H2) holds since
$$\cF G(t,\cdot)(\xi)=\exp\left(-\frac{t|\xi|^2}{2} \right).$$
}
\end{example}

\section{Malliavin differentiability of the solution}
\label{section-Mall-sol}

In this section, we show that the solution of equation \eqref{spde} is Malliavin differentiable and its Malliavin derivative satisfies a certain integral equation. For this, we assume that the function $\sigma$ is affine.

Our first result shows that the sequence of Picard iterations is Malliavin differentiable with respect to $\widehat{N}$ and the corresponding sequence of Malliavin derivatives is uniformly bounded in $L^2(\Omega;\cH)$.

\begin{lemma}
\label{derivative-Picard}
Assume that $\sigma$ is an arbitrary Lipschitz function. Let $(u_n)_{n \geq 0}$ be the sequence of Picard iterations defined by \eqref{def-Picard}. Then  $u_n(t,x) \in \bD^{1,2}$ for any $(t,x) \in [0,T] \times \bR$ and $n \geq 0$, and
$$A:=\sup_{n \geq 0}\sup_{(t,x) \in [0,T] \times \bR}E\|D u_n(t,x)\|_{\cH}^{2}<\infty.$$
\end{lemma}

\noindent {\bf Proof:} {\em Step 1.} We prove that the following property holds for any $n\geq 0$:
\begin{equation}
\label{propertyQ}
\tag{Q}
\left\{
\begin{array}{rcl}
& & \displaystyle u_n(t,x)\in \bD^{1,2} \ \mbox{for any} \ (t,x) \in [0,T] \times \bR,  \quad \mbox{and} \\[2ex]
& & \displaystyle A_n:=\sup_{(t,x) \in [0,T] \times \bR}E\|Du_n(t,x)\|_{\cH}^{2}<\infty
\end{array} \right.
\end{equation}
For this, we use an induction argument on $n$. Property (Q) is clear for $n=0$. We assume that it holds for $n$ and we prove that it holds for $n+1$.

By the definition of $u_{n+1}$ and the fact that the It\^o integral coincides with the Skorohod integral if the integrand is predictable, it follows that
$$u_{n+1}(t,x)=w(t,x)+\int_0^t \int_{\bR}G(t-s,x-y)\sigma(u_n(s,y))L(\delta s,\delta y).$$

We fix $(t,x) \in [0,T] \times \bR$. We apply the Fundamental Theorem of Calculus for the Skorohod integral with respect to $L$ (Theorem \ref{FTC}) to the process:
$$Y(s,y)=G(t-s,x-y)\sigma(u_n(s,y))1_{[0,t]}(s).$$ We need to check that $Y$ satisfies the hypotheses of this theorem.
To check that $Y$ satisfies (i), we apply the Chain Rule (Theorem \ref{chain-rule}) to $F=u_n(s,y)$ and $g=\sigma$. Note that for any $(s,y) \in [0,T] \times \bR$,
\begin{equation}
\label{bound-un}
E|\sigma(u_n(s,y))|^2 \leq 2 D_{\sigma}^2 (1+E|u_n(s,y)|^2)\leq 2 D_{\sigma}^2(1+K_n)<\infty
\end{equation}
 and
$$E\int_0^T \int_{\bR}\int_{\bR_0}|\sigma(u_n(s,y)+D_{r,\xi,z}u_n(s,y))-
\sigma(u_n(s,y))|^2
\nu(dz)d\xi dr \leq$$
$$C_{\sigma}^2 E\int_0^T \int_{\bR}\int_{\bR_0}|D_{r,\xi,z}u_n(s,y)|^2 \nu(dz)d\xi dr \leq C_{\sigma}^2 A_n<\infty,$$
by the induction hypothesis. We conclude that $Y(s,y) \in \bD^{1,2}$ and
\begin{equation}
\label{deriv-Y}
D_{r,\xi,z} Y(s,y)=G(t-s,x-y)\left[\sigma(u_n(s,y)+D_{r,\xi,z}u_n(s,y))-\sigma(u_n(s,y))\right]1_{[0,t]}(s).
\end{equation}

We note that $Y$ satisfies hypothesis (ii) since by \eqref{bound-un},
\begin{eqnarray*}
E\int_0^T \int_{\bR}|Y(s,y)|^2 dyds
&\leq & 2 D_{\sigma}^2 (1+K_n)\int_0^t \int_{\bR}G^2(t-s,x-y)dyds<\infty.
\end{eqnarray*}

To check that $Y$ satisfies hypothesis (iii) i.e. the process $\{D_{r,\xi,z}Y(s,y);s \in [0,T],y \in \bR\}$ is Skorohod integrable with respect to $L$ for any $(r,\xi,z)\in U$, it suffices to show that this process is It\^o intgrable with respect to $L$. Note that $D_{r,\xi,z}u_n(s,y)=0$ if $r>s$ and it is $\cF_s$-measurable if $r \leq s$. Hence, the process
$\{D_{r,\xi,z}Y(s,y);s \in [0,T],y \in \bR\}$ is predictable. By \eqref{deriv-Y} and \eqref{Lip1},
$$E \int_0^T \int_{\bR}|D_{r,\xi,z}Y(s,y)|^2 dyds\leq C_{\sigma}^2 E\int_0^t \int_{\bR} G^2(t-s,x-y)|D_{r,\xi,z}u_n(s,y)|^2 dyds,$$
and hence, %by the induction hypothesis
\begin{eqnarray}
\nonumber
\lefteqn{\int_0^T \int_{\bR}\int_{\bR_0}\left(E \int_0^T \int_{\bR}|D_{r,\xi,z}Y(s,y)|^2 dyds \right)\nu(dz)d\xi dr} \\
\nonumber
& \leq & C_{\sigma}^2 \int_0^t \int_{\bR}G^2(t-s,x-y)E\|D u_n(s,y)\|_{\cH}^2dyds \\
\label{bound-D}
& \leq & C_{\sigma}^2 A_n \int_0^t \int_{\bR}G^2(t-s,x-y)dyds<\infty.
\end{eqnarray}
This proves that $E \int_0^T \int_{\bR}|D_{r,\xi,z}Y(s,y)|^2 dyds<\infty$ for almost all $(r,\xi,z)\in [0,T] \times \bR \times \bR_0$. By Theorem \ref{Ito-Skorohod}.b),
%for almost all $(r,\xi,z)$,
$\{D_{r,\xi,z}Y(s,y);s \in [0,T],y \in \bR\}$ is Skorohod integrable with respect to $L$ and
%its Skorohod integral coincides with its It\^o integral, i.e.
\begin{equation}
\label{Ito-Skorohod-Y-L}
\int_0^T \int_{\bR}D_{r,\xi,z}Y(s,y)L(\delta s,\delta y)=\int_0^T \int_{\bR}D_{r,\xi,z}Y(s,y)L(ds,dy).
\end{equation}

Finally, $Y$ satisfies hypothesis (iv) since by \eqref{Ito-Skorohod-Y-L}, the isometry property \eqref{isometry} and \eqref{bound-D}, we have
\begin{eqnarray*}
\lefteqn{E\int_0^T \int_{\bR} \int_{\bR_0} \left| \int_0^T \int_{\bR} \int_{\bR_0} D_{r,\xi,z} Y(s,y)L(\delta s,\delta y)\right|^2 \nu(dz) d\xi dr =} \\
& & E\int_0^T \int_{\bR} \int_{\bR_0} \left| \int_0^T \int_{\bR} \int_{\bR_0} D_{r,\xi,z} Y(s,y)L(ds,dy)\right|^2 \nu(dz) d\xi dr= \\
& & v  \int_0^T \int_{\bR} \int_{\bR_0} \left(E \int_0^T \int_{\bR} \int_{\bR_0} |D_{r,\xi,z} Y(s,y)|^2 dy ds  \right)\nu(dz) d\xi dr<\infty.
\end{eqnarray*}

By Theorem \ref{FTC}, we infer that $Y \in {\rm Dom}(\delta^L)$, $\delta^L(Y) \in \bD^{1,2}$ and
\begin{equation}
\label{deriv-delta-Y}
D_{r,\xi,z}(\delta^L(Y))=Y(r,\xi)z+\int_0^t \int_{\bR}D_{r,\xi,z}Y(s,y) L(\delta s,\delta y).
\end{equation}
Since $u_{n+1}(t,x)=w(t,x)+\delta^{L}(Y)$, this means that $u_{n+1}(t,x) \in \bD^{1,2}$. Using \eqref{Ito-Skorohod-Y-L} and \eqref{deriv-Y}, we can re-write relation \eqref{deriv-delta-Y} as follows:
\begin{eqnarray}
\nonumber
D_{r,\xi,z} u_{n+1}(t,x)=G(t-r,x-\xi)\sigma(u_n(r,\xi))z+ \int_0^t \int_{\bR}
G(t-s,x-y) \\
\label{deriv-u-n+1}
[\sigma(u_n(s,y)+D_{r,\xi,z} u_n(s,y)) -\sigma(u_n(s,y))]L(ds,dy).
\end{eqnarray}

It remains to prove that
\begin{equation}
\label{An+1finite}
A_{n+1}=\sup_{(t,x) \in [0,T] \times \bR}E\|D u_{n+1}(t,x)\|_{\cH}^2<\infty.
\end{equation}
Using \eqref{deriv-u-n+1}, the isometry property \eqref{isometry}, relation \eqref{bound-un}, and the fact that $\sigma$ is Lipschitz, we see that
\begin{eqnarray*}
\lefteqn{E|D_{r,\xi,z} u_{n+1}(t,x)|^2 \leq 2 z^2 G^2(t-r,x-\xi) E|\sigma(u_n(r,\xi))|^2 } \\
& & +2v E\int_0^t \int_{\bR}G^2(t-s,x-y)|\sigma(u_n(s,y)+D_{r,\xi,z} u_n(s,y)) -\sigma(u_n(s,y))|^2 dy ds \\
& & \leq 4 z^2 D_{\sigma}^2 (1+K_n)G^2(t-r,x-\xi) +2 v C_{\sigma}^2 E\int_0^t \int_{\bR}G^2(t-s,x-y)|D_{r,\xi,z} u_n(s,y)|^2 dy ds.
\end{eqnarray*}
We integrate with respect to $dr d\xi \nu(dz)$ on $[0,T] \times \bR \times \bR_0$. We denote
\begin{equation}
\label{def-nu-t}
\nu_t=\int_0^t \int_{\bR}G^2(s,y)dy ds.
\end{equation}
We obtain:
\begin{eqnarray}
\nonumber
\lefteqn{E\|D u_{n+1}(t,x)\|_{\cH}^2  \leq 4 v D_{\sigma}^2(1+K_n) \nu_t + } \\
\label{bound-D-n+1}
& & 2v  C_{\sigma}^2 \int_0^t \int_{\bR}
G^2(t-s,x-y) E\|D u_n(s,y)\|_{\cH}^2 dy ds\\
\nonumber
& \leq & 4 v D_{\sigma}^2(1+K_n) \nu_t + 2 v C_{\sigma}^2  A_n \nu_t.
\end{eqnarray}
Relation \eqref{An+1finite} follows taking the supremum over $(t,x) \in [0,T] \times \bR$.

{\em Step 2.} We prove that $\sup_{n \geq 1}A_n<\infty$. By
\eqref{bound-D-n+1}, we have:
$$E\|D u_{n+1}(t,x)\|_{\cH}^2 \leq 4 v D_{\sigma}^2(1+K_n) \nu_t + 2vC_{\sigma}^2 \int_0^t H_n(s)J(t-s)ds,$$
where
$V_n(t)=\sup_{x \in \bR} E \|D u_n(t,x)\|_{\cH}^2$ and $J(t)$ is given by \eqref{def-J}. This shows that
$$V_{n+1}(t) \leq 4vD_{\sigma}^2 \nu_T(1+K)+ 2 v C_{\sigma}^2\int_{0}^t V_n(s)J(t-s)ds,$$
where $K$ is given by \eqref{def-K}. %Note that $V_0(t)=0$ for any $t \in [0,T]$.
By Lemma 15 of \cite{dalang99},  $\sup_{n \geq 1}\sup_{t \in [0,T]}V_n(t)<\infty$. $\Box$

\vspace{3mm}
We are now ready to state the main result of the present article.

\begin{theorem}
\label{main}
Assume that $\sigma$ is an affine function, i.e. $\sigma(x)=ax+b$ for some $a,b \in \bR$. If $u$ is the solution of equation \eqref{spde}, then for any $t \in [0,T]$ and $x \in \bR$,
$$u(t,x) \in \bD^{1,2}$$
and the following relation holds in $L^2(\Omega;\cH)$:
\begin{eqnarray}
\nonumber
D_{r,\xi,z} u(t,x)=G(t-r,x-\xi)\sigma(u(r,\xi))z+ \int_0^t \int_{\bR}
G(t-s,x-y) \\
\label{deriv-u}
[\sigma(u(s,y)+D_{r,\xi,z} u(s,y)) -\sigma(u(s,y))]L(ds,dy).
\end{eqnarray}
\end{theorem}

\noindent {\bf Proof:} We fix $(t,x) \in [0,T] \times \bR$.
To prove that $u(t,x) \in \bD^{1,2}$, we apply Theorem \ref{closability-D} to the variables $F_n=u_n(t,x)$ and $F=u(t,x)$. By \eqref{un-conv-u}, $u_n(t,x) \to u(t,x)$ in $L^2(\Omega)$. By Lemma \ref{derivative-Picard}, $u_n(t,x) \in \bD^{1,2}$ for any $n \geq 1$. It remains to prove that $\{Du_n(t,x)\}_{n \geq 1}$ converges in $L^2(\Omega;\cH)$. Let
$$M_n(t)=\sup_{x \in \bR}E\|Du_n(t,x)-Du_{n-1}(t,x)\|_{\cH}^2.$$

We write relation \eqref{deriv-u-n+1} for $D_{r,\xi,z}u_{n+1}(t,x)$ and $D_{r,\xi,z}u_n(t,x)$. We take the difference between these two equations. We obtain:
\begin{eqnarray}
\nonumber
\lefteqn{D_{r,\xi,z}u_{n+1}(t,x)-D_{r,\xi,z}u_n(t,x)= G(t-r,x-\xi)[\sigma(u_n(r,\xi))-\sigma(u_{n-1}(r,\xi))]z+ }\\
\nonumber
& & \int_0^t \int_{\bR} G(t-s,x-y) \{[\sigma(u_n(s,y)+D_{r,\xi,z}u_n(s,y))-\sigma(u_n(s,y))]-\\
\label{Malliavin}
& & \quad \quad \quad [\sigma(u_{n-1}(s,y)+D_{r,\xi,z}u_{n-1}(s,y))-\sigma(u_{n-1}(s,y))] \} L(ds,dy).
\end{eqnarray}

At this point, we use the assumption that $\sigma$ is the affine function $\sigma(x)=ax+b$. (An explanation why this argument does not work in the general case is given in Remark \ref{rem-sigma-affine} below.) In this case, relation \eqref{Malliavin} has the following simplified expression:
\begin{eqnarray*}
\nonumber
\lefteqn{D_{r,\xi,z}u_{n+1}(t,x)-D_{r,\xi,z}u_n(t,x)= a G(t-r,x-\xi)[u_n(r,\xi)-u_{n-1}(r,\xi)]z+ }\\
& & a \int_0^t \int_{\bR}G(t-s,x-y)[D_{r,\xi,z}u_n(s,y)-D_{r,\xi,z}u_{n-1}(s,y)]L(ds,dy).
\end{eqnarray*}
Using It\^o's isometry and the inequality $(a+b)^2 \leq 2(a^2+b^2)$, we obtain:
\begin{eqnarray*}
\lefteqn{E|D_{r,\xi,z}u_{n+1}(t,x)-D_{r,\xi,z}u_n(t,x)|^2 \leq 2a^2 z^2 G^2(t-r,x-\xi)b_n^2+} \\
& & 2a^2 v E\int_0^t \int_{\bR}G^2(t-s,x-y)|D_{r,\xi,z}u_n(s,y)-D_{r,\xi,z}u_{n-1}(s,y)|^2 dyds,
\end{eqnarray*}
where
$b_n^2=\sup_{(s,y) \in [0,T] \times \bR}E|u_n(s,y)-u_{n-1}(s,y)|^2$. Note that both sides of the previous inequality are zero if $r>t$. Taking the integral with respect to $dr d\xi \nu(dz)$ on $[0,T] \times \bR \times \bR_0$,
we obtain:
\begin{eqnarray*}
\lefteqn{E\|Du_{n+1}(t,x)-Du_n(t,x)\|_{\cH}^2 \leq 2a^2 v  \nu_t b_n^2+} \\
& & 2a^2 v E\int_0^t \int_{\bR}G^2(t-s,x-y)E\|Du_n(s,y)-Du_{n-1}(s,y)\|_{\cH}^2 dyds,
\end{eqnarray*}
where $\nu_t$ is given by \eqref{def-nu-t}. Recalling the definition of
$M_n(t)$, we infer that
$$M_{n+1}(t) \leq C_n+2a^2 v \int_0^t M_n(s)J(t-s)ds,$$
where $C_n=2 a^2 v \nu_t b_n^2 $ and the function $J$ is given by \eqref{def-J}.
By relation \eqref{sum-converges}, we know that $\sum_{n \geq 1}b_n<\infty$, which means that $\sum_{n \geq 1} C_n^{1/2}<\infty$. By Lemma \ref{Dalang-lemma} (Appendix A), we conclude that
$$\sum_{n \geq 1}\sup_{t \leq T}M_n(t)^{1/2}<\infty.$$
Hence, the sequence $\{Du_n(t,x)\}_{n \geq 1}$ converges in $L^2(\Omega;\cH)$ to a variable $U(t,x)$, uniformly in $(t,x) \in [0,T] \times \bR$. By Theorem
\ref{closability-D}, $u(t,x) \in \bD^{1,2}$ and $Du_n(t,x) \to D u(t,x)$ in $L^2(\Omega;\cH)$. Hence,
$$\sup_{(t,x) \in [0,T] \times \bR} E\|Du_n(t,x)-D u(t,x)\|_{\cH}^2 \to 0.$$
Relation \eqref{deriv-u} follows by taking the limit in $L^2(\Omega;\cH)$ as $n \to \infty$ in \eqref{deriv-u-n+1}. $\Box$

\begin{remark}
\label{rem-sigma-affine}
{\rm Unfortunately, we were not able to extend Theorem \ref{main} to an arbitrary Lipschitz function $\sigma$. To see where the difficulty comes from, recall that we need to prove that $\{Du_n(t,x)\}_{n \geq 1}$ converges in $L^2(\Omega;\cH)$, and the difference $D_{r,\xi,z}u_{n+1}(t,x)-D_{r,\xi,z}u_n(t,x)$ is given by \eqref{Malliavin}. For an arbitrary Lipschitz function $\sigma$, by relation \eqref{Lip1}, we have:
\begin{eqnarray*}
\lefteqn{|\sigma(u_n(s,y)+D_{r,\xi,z}u_n(s,y))-\sigma(u_{n-1}(s,y)+D_{r,\xi,z}u_{n-1}(s,y))|
\leq } \\
& & C_{\sigma}|(u_n(s,y)+D_{r,\xi,z}u_n(s,y))-(u_{n-1}(s,y)+D_{r,\xi,z}u_{n-1}(s,y))| \leq \\
& & C_{\sigma}\{|u_n(s,y)-u_{n-1}(s,y)|+|D_{r,\xi,z}u_n(s,y))-D_{r,\xi,z}u_{n-1}(s,y)|\}.
\end{eqnarray*}

Using \eqref{Malliavin}, the isometry property \eqref{isometry}, the inequality $(a+b)^2 \leq 2(a^2+b^2)$, and the previous inequality, we have:
\begin{eqnarray*}
\lefteqn{E|D_{r,\xi,z}u_{n+1}(t,x)-D_{r,\xi,z}u_n(t,x)|^2 \leq 2
z^2 C_{\sigma}^2 G^2(t-r,x-\xi)E|u_n(r,\xi)-u_{n-1}(r,\xi)|^2  } \\
& & +4 v C_{\sigma}^2  \int_0^t \int_{\bR}G^2(t-s,x-y)E|u_n(s,y)-u_{n-1}(s,y)|^2 dy ds\\
& & +4 v C_{\sigma}^2 \int_0^t \int_{\bR}G^2(t-s,x-y) E|D_{r,\xi,z}u_n(s,y))-D_{r,\xi,z}u_{n-1}(s,y)|dy ds.
\end{eqnarray*}
The problem is that the second term on the right-hand side of the inequality above does not depend on $(r,\xi,z)$ and hence its integral with respect to $dr d\xi \nu(dz)$ on $[0,T] \times \bR \times \bR_0$ is equal to $\infty$. }
\end{remark}

\appendix

\section{A variant of Gronwall's lemma}
\label{section-app}

The following result is a variant of Lemma 15 of \cite{dalang99}, which is used in the proof of Theorem \ref{main}.

\begin{lemma}
\label{Dalang-lemma}
Let $(f_n)_{n \geq 0}$ be a sequence of non-negative functions defined on $[0,T]$ such that $M=\sup_{t \in [0,T]}f_0(t)<\infty$ and for any $t \in [0,T]$ and $n \geq 0$,
$$f_{n+1}(t) \leq C_n+\int_0^t f_n(s) g(t-s)ds,$$
where $g$ is a non-negative function on $[0,T]$ with $\int_0^T g(t)dt<\infty$ and $(C_n)_{n \geq 0}$ is a sequence of non-negative constants. Then, there exists a sequence $(a_n)_{n \geq 0}$ of non-negative constants which satisfy $\sum_{n \geq 0}a_n^{1/p}<\infty$ for any $p>1$, such that for any $t \in [0,T]$ and $n \geq 0$,
\begin{equation}
\label{Gronwall-rel}
f_n(t) \leq C_n +\sum_{j=1}^{n-1}C_j a_{n-j}+C_0 a_n M.
\end{equation}
In particular, if $\sum_{n \geq 1}C_n^{1/p}<\infty$ for some $p>1$, then
$$\sum_{n \geq 1}\sup_{t \in [0,T]}f_n(t)^{1/p}<\infty.$$
\end{lemma}

\noindent {\bf Proof:} Let $G(T)=\int_0^T g(t)dt$, $(X_i)_{i \geq 1}$ be a sequence of i.i.d. random variables on $[0,T]$ with density function $g(t)/G(T)$, and $S_n=\sum_{i=1}^n X_i$. Following exactly the same argument as in the proof of Lemma 15 of \cite{dalang99}, we have:
\begin{eqnarray*}
f_n(t) & \leq & C_n+C_{n-1}G(T)P(S_1 \leq t)+\ldots+C_1 G(T)^{n-1}P(S_{n-1}\leq t) +\\
& & C_0G(T)E[1_{\{S_n \leq t\}}f_0(t-S_n)].
\end{eqnarray*}
Relation \eqref{Gronwall-rel} follows with $a_n=G(T)^nP(S_n \leq T)$ for $n \geq 1$.
The fact that $\sum_{n \geq 1}a_n^{1/p}<\infty$ for all $p \geq 1$ was shown in the proof of Lemma 15 of \cite{dalang99}.

To prove the last statement, we let $a_0=1$ and $M_1=\max(M,1)$. Then $f_n(t) \leq M_1 \sum_{j=0}^n C_j a_{n-j}$ and hence, $\sup_{t \leq T}f_n(t)^{1/p} \leq M_1^p \sum_{j=0}^{n}C_j^{1/p}a_{n-j}^{1/p}$. We conclude that
$$\sum_{k=0}^{n}\sup_{t \leq T}f_k(t)^{1/p} \leq M_1^p \sum_{j=0}^{n}C_j^{1/p}\sum_{k=j}^{n}a_{k-j}^{1/p}\leq M_1^p \sum_{j\geq 0}C_j^{1/p}\sum_{k \geq 0}a_{k}^{1/p}:=C<\infty.$$
$\Box$

\end{document}